\documentclass[a4paper,UKenglish,cleveref, autoref, thm-restate]{lipics-v2021}

\usepackage{algorithm}
\usepackage{algpseudocode}
\algrenewcommand\algorithmicrequire{\textbf{Input:}}
\algrenewcommand\algorithmicensure{\textbf{Output:}}

\usepackage{placeins}
\usepackage{booktabs}
\usepackage{accents}
\newcommand{\ubar}[1]{\underaccent{\bar}{#1}}

\newcommand{\tabtop}{\rule{0pt}{9pt}}

\title{Identifying Multi-Hit Cancer Drivers Without Massive Parallelization: A CP, MIP, and Column Generation Framework} %TODO Please add

\titlerunning{Identifying Multi-Hit Cancer Drivers Without Massive Parallelization} %TODO optional, please use if title is longer than one line

\author{Rick S. H. Willemsen}{Singapore University of Technology and Design, Engineering Systems and Design, Singapore}{rick_willemsen@sutd.edu.sg}{https://orcid.org/0000-0002-7275-4923}{}%TODO mandatory, please use full name; only 1 author per \author macro; first two parameters are mandatory, other parameters can be empty. Please provide at least the name of the affiliation and the country. The full address is optional. Use additional curly braces to indicate the correct name splitting when the last name consists of multiple name parts.

\author{Tenindra Abeywickrama}{RIKEN Center for Computational Science, Japan}{tenindra.abeywickrama@riken.jp}{https://orcid.org/0000-0003-1053-2193}{}

\author{Ramu Anandakrishnan}{Virginia Tech, Edward Via College of Osteopathic Medicine, USA}{ramu@vt.edu}{https://orcid.org/0000-0003-0422-3984}{}

\authorrunning{R.\,S.\,H. Willemsen, T. Abeywickrama, and R. Anandakrishnan} %TODO mandatory. First: Use abbreviated first/middle names. Second (only in severe cases): Use first author plus 'et al.'

\Copyright{Rick S. H. Willemsen, Tenindra Abeywickrama, and Ramu Anandakrishnan} %TODO mandatory, please use full first names. LIPIcs license is "CC-BY";  http://creativecommons.org/licenses/by/3.0/

 %TODO mandatory: Please choose ACM 2012 classifications from https://dl.acm.org/ccs/ccs_flat.cfm 

\relatedversion{} %optional, e.g. full version hosted on arXiv, HAL, or other respository/website
\supplement{}
\supplementdetails[linktext={}, subcategory={Source code}]{Software}{https://github.com/rickwillemsen/MHCDSCP}
\acknowledgements{We thank the anonymous reviewers for their helpful comments that significantly improved our work and the authors of BiGPICC~\cite{oles2025bigpicc} for sharing the datasets.}%optional

\EventEditors{John Q. Open and Joan R. Access}
\EventNoEds{2}
\EventLongTitle{42nd Conference on Very Important Topics (CVIT 2016)}
\EventShortTitle{CVIT 2016}
\EventAcronym{CVIT}
\EventYear{2016}
\EventDate{December 24--27, 2016}
\EventLocation{Little Whinging, United Kingdom}
\EventLogo{}
\SeriesVolume{42}
\ArticleNo{23}
\begin{document}

\nolinenumbers
\maketitle

\begin{abstract}
Cancer is often driven by specific combinations of an estimated two to nine gene mutations, known as multi-hit combinations. Identifying these multi-hit combinations of gene mutations that drive cancer is critical for understanding carcinogenesis and designing targeted therapies. We formalize this challenge as the Multi-Hit Cancer Driver Set Cover Problem (MHCDSCP), optimizing the selection of gene combinations to maximize tumor coverage while strictly minimizing normal sample misclassification. While existing approaches rely on exhaustive enumeration and massive parallelization, we introduce fast heuristics based on constraint programming and mixed integer programming formulations. Evaluated on real-world cancer genomics data, our framework matches state-of-the-art supercomputing methods using a single commodity CPU in under a minute. We also propose a price-and-branch heuristic which, by solving the root node to optimality, provides the first provably optimal solutions for over half of the benchmark instances, thereby verifying the near-optimality of our fast heuristics. These findings demonstrate that on real-world problem instances, the MHCDSCP is far less computationally demanding than previously believed, providing an accessible baseline that enables the exploration of previously intractable multi-hit modeling assumptions.
\end{abstract}

\section{Introduction}
Cancer remains a leading cause of mortality worldwide, with an estimated 9.7 million cancer-related deaths in 2022 \cite{bray2024global}. This number is projected to rise to over 18 million in 2050 \cite{luo2025global}. Understanding the biological mechanisms that drive cancer development is therefore essential. According to the multi-hit theory of carcinogenesis, multiple mutated genes (called hits) are required to cause cancer in humans \cite{armitage1954age, nordling1953new, sutherl1984multihit}. Mathematical models estimate that the number of mutated genes (hit size) typically ranges from two to nine, depending on the type of cancer \cite{little2003stochastic, zhang2005estimating, tomasetti2015only, ashley1969two, anandakrishnan2019estimating}. Despite this understanding, identifying the specific combinations of gene mutations responsible for a particular cancer remains a challenge. Discovering these gene combinations leads to a better understanding of carcinogenesis and enables the development of more effective targeted therapies \cite{al2012combinatorial, dash2019differentiating}. Although other factors, such as the tumor environment \cite{schneider2017tissue} and epigenetic modifications \cite{stahl2016epigenetics}, also influence cancer development, this work focuses on identifying carcinogenic combinations of mutated genes.

A popular model for identifying carcinogenic multi-hit gene combinations involves formulating the task as a weighted set cover (WSC) problem, as introduced by \cite{dash2019differentiating}. In this model, the goal is to find a set of gene combinations that collectively ``explain'' mutations in tumor samples, while avoiding mutations in normal samples. By penalizing the inclusion of normal samples, the model isolates potential cancer driver genes from passenger mutations. Specifically, the effectiveness of a single combination is measured by the weighted difference between the number of tumor and normal samples covered (penalizing the latter). The resulting set of gene combinations is then evaluated as a binary classifier. While the task is framed as a WSC problem, the existing literature \cite{dash2019differentiating, al2020identifying, dash2021scaling, oles2025bigpicc, prabhu2026looking}, does not optimize a global objective function, but instead employs an iterative greedy heuristic. In each iteration, the algorithm selects a gene combination with the highest current weight. All tumor samples covered by the selected gene combination are then removed and the weights of all remaining combinations are updated. This process continues until all tumor samples are covered. However, we identify two key limitations in current work adopting this model:

\begin{enumerate}

\item To the best of our knowledge, all prior work that adopted this model relies on massive parallelization using supercomputing infrastructure. This is primarily due to the exhaustive enumeration of all or large subsets of possible gene combinations. For instance, one approach employed 6,000 V100 high-performance data center GPUs (reported as 48 million GPU cores) \cite{dash2021scaling}, while the current state-of-the-art method used 280 nodes on the Summit supercomputer (estimated 11,760 CPU cores) \cite{oles2025bigpicc}. A concurrent preprint \cite{prabhu2026looking} reported using 3072 nodes (147,456 CPU cores) on the Fugaku supercomputer.

\item These approaches adopt a greedy set cover heuristic that enforces strict coverage of tumor samples, which presents several limitations. First, strict coverage risks overfitting to noisy or misclassified tumor samples. Second, the greedy nature of these algorithms does not allow flexibility in other solution quality measures, such as minimizing the number of selected gene combinations. Third, the strict coverage assumption is misaligned with some of the evaluation criteria used. The symmetric evaluation criteria used in previous works, give equal weight to misclassified tumor and normal samples, while the greedy approach prioritizes covering all tumor samples regardless of classification accuracy.
\end{enumerate}

To overcome these limitations, and given that resulting multi-hit gene combinations are ultimately evaluated as a binary classifier, we focus on directly optimizing the classification accuracy. A subset of gene combinations induces a binary classification of samples: a sample is classified as a tumor if at least one selected gene combination occurs in it, or classified as a normal sample otherwise. Therefore, we propose an objective function that balances two criteria, namely maximizing the number of correctly classified tumor samples (true positives) and minimizing the number of times normal samples are misclassified as tumor (false positives). We formalize this problem as the Multi-Hit Cancer Driver Set Cover Problem (MHCDSCP). Note that for both the WSC and the MHCDSCP it is not necessary to enumerate all possible gene combinations. Instead, the main goal is to identify which subset of gene combinations leads to the best classification performance. Our contributions are as follows:

\begin{enumerate}

\item To model the MHCDSCP, we propose both a constraint programming (CP) and a mixed integer programming (MIP) formulation, which take as input a set of candidate gene combinations. When evaluated on real-world cancer genomics data, both models find combinations that match the classification performance of existing massively parallel approaches while requiring only a single commodity CPU. Specifically, the MIP formulation runs within 1 minute and the CP model in 20 minutes.

\item We also introduce a price-and-branch heuristic based on column generation that dynamically generates gene combinations as needed. Within a 5-minute time limit, this approach solves over half of the benchmark instances to provable optimality for the first time and demonstrates that our MIP formulation yields near-optimal solutions. Notably, our work is the first to provide meaningful insights into the solution quality of identified gene combinations, in stark contrast to the lack of optimality guarantees in previous greedy WSC approaches.

\item Our framework enables flexible optimization goals and constraints, including optimizing the number of combinations and accommodating flexibility in covering tumors. By relaxing the assumption that all tumor samples must be covered, our approach can mitigate the influence of noisy or misclassified samples. Conversely, strict coverage can be easily enforced if required. Furthermore, our methods achieve similar classification performance to existing methods while identifying smaller sets of gene combinations, potentially reaffirming the multi-hit theory that a limited number of gene mutations drive the development of cancer.

\end{enumerate}

Our key finding, that near-optimal solutions for real-world problem instances can be obtained without use of massive parallelization, indicates that the MHCDSCP, and its variations, are not as computationally demanding as previously believed. Moreover, our work provides an easily accessible baseline for researchers working in cancer genomics. Crucially, our models enable the evaluation of hypotheses in multi-hit cancer theory that were previously intractable and identify potential directions for cancer biology research, such as determining the biologically aligned number of combinations. In particular, the CP formulation provides a platform for exploration of complex biological constraints that may be difficult to linearize. 

The remainder of this paper is organized as follows. In \Cref{section:relatedwork} we further discuss related work. In \Cref{section2:problemdescription} we introduce the MHCDSCP and provide its mathematical formulations. \Cref{section3:colgen} details how these formulations can be solved within a column generation-based framework. The proposed methods are evaluated in \Cref{section4:results} on real-world data. Finally, \Cref{section5:conclusion} concludes and discusses directions for further research.

\section{Further Related Work}\label{section:relatedwork}

\textbf{Alternatives to Multi-Hit Models:} While several approaches aim to identify individual driver mutations \cite{tian2014contrastrank, tamborero2013oncodriveclust, dees2012music, kumar2016unsupervised}, such gene mutations alone are generally insufficient to explain cancer development, of which several examples are provided by \cite{dash2019differentiating}. Another theory is based on the assumption of mutual exclusivity, which states that mutations in certain genes rarely co-occur and that a mutation in any of them may be sufficient to drive cancer. Mutual exclusivity is identified through so-called driver pathway discovery \cite{vandin2012novo}. Based on this assumption multiple methods have been introduced \cite{leiserson2013simultaneous, leiserson2015comet, park2022superdendrix, kim2017wesme, dao2017bewith, zhang2017discovery}. In particular, mixed integer programming has been successfully applied to this problem \cite{zhao2012efficient}. Other methods focus on identifying modules of co-occurring cancer driver genes using biological interaction networks \cite{klein2021identifying, zhang2024identifying, zhang2024identifying}. All aforementioned approaches differ in the modeling assumptions from multi-hit models, which aim to identify multi-hit gene combinations.

\textbf{Enumeration-Based Weighted Set Cover Approaches:} In the existing literature, the weighted set cover problem formulation for discovering multi-hit carcinogenic gene combinations is primarily addressed using a two-step process that enumerates combinations and greedily selects the combination with maximum weight. Early approaches involved exhaustively enumerating all gene combinations for 2, 3 and 4 hits \cite{dash2019differentiating, al2020identifying, dash2021scaling}. BiGPICC ~\cite{oles2025bigpicc}, a graph-based heuristic, was able to identify gene combinations up to 8 hits by limiting exhaustive enumeration to a large subset of candidate combinations chosen by a clustering process. Subsequent to our initial submission, a concurrent preprint \cite{prabhu2026looking} revisited enumerating all gene combinations but with pruning of combinations that cannot cover any tumors. However, this approach does not support ranges of hits, requires substantial computational resources (nearly 150,000 cores) to evaluate up to 4-hits, and was not compared to BiGPICC.

\textbf{Set Cover Algorithms:} Constraint programming and mixed integer programming tools are well-known for solving set cover-type problems \cite{balas1972set, hochbaum1997approximating, hnich2006constraint}, and have been widely used in constraint-based classification and pattern mining tasks \cite{guns2011itemset, de2010constraint}. In particular, optimization-based techniques have been successfully applied to bioinformatics problems \cite{han2007frequent}. While machine learning is highly effective in building classifiers for cancer genomics~\cite{eraslan2019deep}, set cover is still preferred as our objective is to identify explicit multi-hit gene combinations that are then evaluated as classifiers. This cannot be achieved by black-box ML models that lack interpretability. Furthermore, recent ML-based set cover approaches such as Graph-SCP~\cite{shafi2025graphscp}, while promising, are designed to augment rather than replace traditional solvers.

\section{Problem Description}\label{section2:problemdescription}
In this section, we introduce the Multi-Hit Cancer Driver Set Cover Problem (MHCDSCP) as a binary classification problem. Afterwards, we provide a constraint programming (CP) and mixed integer programming (MIP) formulation for the MHCDSCP.

\subsection{Notation}
\begin{figure}[h]
    \centering
    \includegraphics[width=0.8\linewidth]{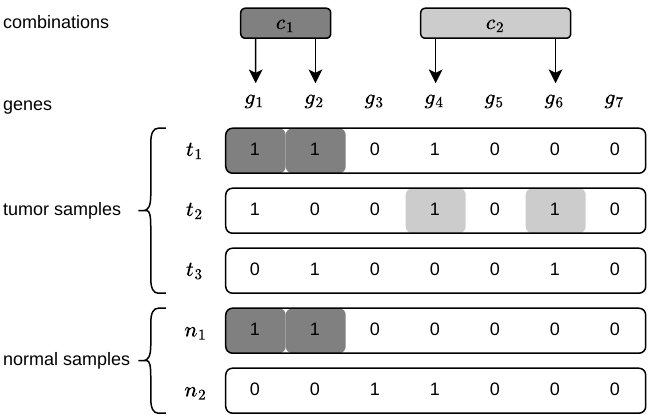}
    \caption{An example with seven genes and five samples represented as a binary matrix. Two gene combinations with a hit size of two are selected, namely $c_1$ in dark grey and $c_2$ in light grey. Together, they cover tumor samples $t_1$ and $t_2$ (true positives) and normal sample $n_1$ (false positive).}
    \label{fig:problemDescription}
\end{figure}

Let $G$ be the set of genes and $S$ the set of samples, partitioned into tumor samples $T$ and normal samples $N$, with $S=T \cup N$. The input data are represented by a binary matrix $\mathbf{M}$ of size $|S|\times |G|$, where entry $\mathbf{M}_{sg}=1$ indicates that gene $g\in G$ is mutated in sample $s\in S$, and $\mathbf{M}_{sg}=0$ otherwise (as in, e.g., \cite{oles2025bigpicc}). This binary matrix is visualized in \Cref{fig:problemDescription}.

A gene combination is a subset of genes whose cardinalities lie within a predefined hit range. Let $\ubar{k}$ and $\bar{k}$ denote the minimum and maximum hit sizes, respectively. We define $C$ as the set of all possible gene combinations such that the cardinality lies within the given hit range. A gene combination $c\in C$ covers a sample $s\in S$ if all genes in combination $c$ are mutated in $s$. Furthermore, we denote by $C_s\subseteq C$ all the combinations that cover sample $s$ and we denote by $S_c\subseteq S$ all the samples that are covered by combination $c$.

\subsection{Multi-Hit Cancer Driver Set Cover Problem}
In the MHCDSCP, we aim to select a subset of gene combinations $C^*\subseteq C$ that accurately classifies samples as tumor or normal. Selecting a set $C^*$ induces a binary classification of samples as follows. A sample $s\in S$ is classified as tumor if at least one $c\in C^*$ covers~$s$, and as a normal sample otherwise. The MHCDSCP seeks to identify a subset of gene combinations that balances two criteria, namely maximizing the number of correctly classified tumor samples and minimizing the number of times a normal sample is misclassified.

Let $T^*\subseteq T$ and $N^*\subseteq N$ be the set of samples that are classified by the selected gene combinations $C^*$ as tumor samples. We define the number of true positives as $|T^*|$, false negatives as $|T\setminus T^*|$, false positives as $|N^*|$, and true negatives as $|N\setminus N^*|$. \Cref{fig:problemDescription} illustrates selected gene combinations that cover a number of samples, leading to a number of true positives and false positives.

\subsection{Constraint Programming Formulation}
Let $x_t$ be a binary variable that takes the value 1 if tumor sample $t\in T$ is classified as tumor. Let the integer variable $y_n$ represent the number of times a normal sample $n \in N$ is misclassified as tumor. We define a binary variable $z_c$ that takes value 1 if gene combination $c\in C$ is selected. The maximum number of selected combinations is denoted by the parameter~$\beta$. 

The selection of gene combinations can be formulated using CP. Since typical performance measures (see Appendix \ref{section2.3:performancemeasures}) are nonlinear, we instead define a linear objective that balances two criteria, namely maximizing the number of correctly classified tumor samples (true positives) and minimizing the number of times a normal sample is misclassified (approximating false positives). The CP formulation for the MHCDSCP is as follows:
\begin{align}
\max\ & \sum_{t\in T} x_t -  \sum_{n\in N} y_n, \label{cp:cons0} \\
   \text{s.t.}\  & x_t =\max_{c \in C_t} \{z_c\}, && \forall t \in T, \label{cp:cons1} \\
    & y_n = \sum_{c \in C_n} z_c, && \forall n \in N, \label{cp:cons2} \\
    & \sum_{c \in C} z_c \le \beta, \label{cp:cons3} \\
    & x_t \in \{0,1\}, && \forall t\in T ,\label{cp:cons4}\\
    & y_n \in \{0, 1, \dots, \beta\}, && \forall n\in N ,\label{cp:cons5}\\
    & z_c \in \{0,1\},   && \forall c\in C.\label{cp:cons6}
\end{align}
The objective (\ref{cp:cons0}) maximizes the correctly classified tumor samples while minimizing misclassification of normal samples. The max-equality constraints~(\ref{cp:cons1}) ensure that a tumor sample $t$ is classified as tumor if it is covered by at least one selected gene combination. Constraints~(\ref{cp:cons2}) count the number of times a normal sample $n$ appears in a selected combination. The maximum number of selected gene combinations is enforced by (\ref{cp:cons3}). The remaining constraints (\ref{cp:cons4})-(\ref{cp:cons6}) model the domain of the variables.

\subsection{Mixed Integer Programming Formulation}
Using the same notation as the CP formulation, we also introduce an MIP formulation. The two formulations are nearly identical, with the main difference being that the max-equality constraints are linearized, which determine whether a tumor sample is correctly classified. This linearization allows the problem to be solved using a standard MIP solver. The MIP is formulated as follows:
\begin{align}
    \max\ & \sum_{t\in T} x_t - \sum_{n\in N} y_n ,  \label{form1:cons1} \\
    \text{s.t.}\ & x_t - \sum_{c\in C_t} z_c\leq 0, && \forall t\in T ,\label{form1:cons2} \\
    & \sum_{c\in C_n} z_c-  y_n   =0 , && \forall n\in N, \label{form1:cons3} \\
    & \sum_{c\in C} z_c  \leq  \beta  ,\label{form1:cons4} \\
    & x_t \in  \{0,1\} ,&& \forall t\in T ,\label{form1:cons5} \\
    & y_n \in  \{0, 1, \dots, \beta\},&& \forall n\in N ,\label{form1:cons6} \\
    & z_c \in \{0,1\} , && \forall c\in C .\label{form1:cons7}
\end{align}
Note that (\ref{form1:cons2}) differs from (\ref{cp:cons1}) in the CP formulation. In an optimal solution, both constraints provide the same result, since $x_t$ is maximized in the objective function. Furthermore, we relax the integrality of the $x_t$ and $y_n$ variables when solving the model, because the integrality is implied via the binary $z_c$ variables. Preliminary experiments show that this improves computational performance. 
 
The CP formulation (\ref{cp:cons0})-(\ref{cp:cons6}) and the MIP formulation (\ref{form1:cons1})-(\ref{form1:cons7}) can be solved only when considering a limited set of gene combinations. In general, the set $C$ contains all gene combinations of size up to the maximum hit size $\bar{k}$, and its cardinality increases exponentially with $\bar{k}$. As a result, both formulations may contain an intractably large number of variables. 

Preliminary computational results indicate that MIP is preferred over CP. Column generation is a technique specifically designed to address problems with an exponential number of variables while still guaranteeing optimality. Next, we present a column generation framework to efficiently solve the MIP.

\section{Column Generation-Based Framework}\label{section3:colgen}
The large number of gene combinations $C$ makes it impractical to solve the MIP formulation (\ref{form1:cons1})-(\ref{form1:cons7}) directly. Column generation is a well-established technique for addressing such large-scale problems by working with a restricted subset of variables and iteratively introducing new variables as needed \cite{desaulniers2006column}. We first present a relaxed version of the MIP. We then introduce a fast heuristic and a price-and-branch heuristic to provide bounds on the objective value.

\subsection{Restricted Master Problem}
In column generation, we start from a limited set of gene combinations $\Tilde{C}\subseteq C$, which may initially be empty. Instead of solving the MIP (\ref{form1:cons1})-(\ref{form1:cons7}), we solve a so-called restricted master problem, which corresponds to the linear relaxation of the MIP restricted to gene combinations in $\tilde{C}$.

Solving the restricted master problem provides information in the form of dual values that can be used to generate gene combinations not contained in $\tilde{C}$ that may improve the objective value. Such gene combinations correspond to new columns and are added to $\tilde{C}$ in an iterative procedure. Repeating this process yields an optimal (possibly fractional) solution to the linear relaxation of the MIP (\ref{form1:cons1})-(\ref{form1:cons7}). The resulting objective value provides a valid upper bound on the optimal objective value of the original problem.

We relax the integrality constraints on $z_c$ (\ref{form1:cons6}) to $z_c\geq 0$, noting that $z_c$ will not exceed one in an optimal solution. Furthermore, by relaxing the equality constraints (\ref{form1:cons3}) into inequalities, we obtain the following restricted master problem
\begin{align}
    \max\ & \sum_{t\in T} x_t  - \sum_{n\in N} y_n , \label{form2:cons1} \\
    \text{s.t.}\ & x_t - \sum_{c\in \tilde{C}_t} z_c \leq 0, && \forall t\in T ,\label{form2:cons2} \\
    & \sum_{c\in \tilde{C}_n} z_c -  y_n \leq 0 , && \forall n\in N, \label{form2:cons3} \\
    & \sum_{c\in \tilde{C}} z_c  \leq  \beta  ,\label{form2:cons4} \\
    & x_t \in  [0, 1], && \forall t\in T ,\label{form2:cons5} \\
    & y_n \geq 0 ,&& \forall n\in N ,\label{form2:cons6} \\
    & z_c \geq 0 , && \forall c\in \tilde{C} .\label{form2:cons7} 
\end{align}

\subsection{Greedy Restricted MIP Heuristic}
We first present a Greedy Restricted MIP (GR-MIP) heuristic that makes use of the column generation interpretation to quickly obtain feasible solutions for the MIP (\ref{form1:cons1})-(\ref{form1:cons7}). Note that this heuristic can also solve the CP formulation, in which case we refer to it as GR-CP.

The restricted master problem is initialized by adding a set of gene combinations $\tilde{C}$ generated by a greedy procedure. After adding these columns to the model and restricting the variables $z_c$ to be binary again, the resulting MIP can be solved using a standard MIP solver. Any feasible solution obtained provides a lower bound on the objective value. Recall that due to the binary variables $z_c$ and the direction of the objective function, the variables $x_t$ and $y_n$ are forced to take integer values. 

The greedy restricted procedure for generating gene combinations is described in \Cref{alg:greedyprocedure}. First, we select $\gamma_1$ candidate genes that are most frequently mutated in tumor samples. Second, a hit size $k\in [\ubar{k}, \bar{k} ] \cap \mathbb{Z}$ is chosen randomly. Third, a gene combination of size $k$ is generated by randomly sampling from the selected candidate genes. This procedure is repeated until at most $\gamma_2$ unique gene combinations (columns) have been generated.

\begin{algorithm}
\caption{Greedy generation of gene combinations}
\label{alg:greedyprocedure}
\begin{algorithmic}[1]
\Require Binary matrix $\mathbf{M}$
\Ensure Set of gene combinations

\State Rank genes by mutation frequency in tumor samples
\item Select the top $\gamma_1$ genes as candidate set

\Repeat
    \State Randomly select a hit size $k\in [\ubar{k}, \bar{k} ] \cap \mathbb{Z}$
    \State Randomly select $k$ genes from the candidate set to form a gene combination
\Until{$\gamma_2$ unique gene combinations have been generated}

\end{algorithmic}
\end{algorithm}
\FloatBarrier

\subsection{Price-and-Branch Heuristic}
Instead of adding greedily generated columns, we employ a so-called price-and-branch heuristic \cite{sadykov2019primal}. This is a column generation-based method consisting of two phases. In the first phase, it solves the root node to optimality using column generation to identify promising gene combinations. In the second phase, integrality constraints (of $z_c$ variables) are restored and the resulting problem is solved using a standard MIP solver. A column generation-based method allows us to efficiently explore the set of all gene combinations without the need for full enumeration.

Column generation works as follows. First, the restricted master problem is solved and the dual values associated with the constraints are passed to the pricing problem. These dual values are used in the pricing problem to identify gene combinations (columns) with positive reduced costs. When these positive reduced cost columns are added to $\tilde{C}$, these columns may improve the objective function. 

The processes of solving the restricted master problem and the pricing problem are repeated iteratively until no more gene combinations (columns) with positive reduced cost exist. At each iteration, we obtain a (possibly fractional) solution ($x^*, y^*, z^*$). We apply a simple rounding heuristic to obtain a valid lower bound on the objective value, which is summarized in \Cref{alg:roundingheuristic}. When there are no further positive reduced cost combinations, an optimal solution to the linear relaxation of the MIP is obtained, providing a valid upper bound on the objective value. In some cases, this solution is already integral, yielding a provable optimal solution to the original problem. Otherwise, we add integrality restrictions to $z_c$ (\ref{form2:cons7}) again and solve it using an MIP solver.

\begin{algorithm}
\caption{Rounding heuristic}
\label{alg:roundingheuristic}
\begin{algorithmic}[1]
\Require Fractional solution $x^*, y^*, z^*$
\Ensure Feasible integral solution and corresponding lower bound $LB$

\State $k \gets \lceil \sum_{c\in \tilde{C}} z_c^* \rceil$ 

\State Select the $k$ largest $z_c^*$ to obtain a feasible integral solution

\State Set $LB$ equal to the corresponding objective value

\end{algorithmic}
\end{algorithm}
\FloatBarrier

\subsubsection{Pricing Problem for the Identification of Gene Combinations}
The goal of the pricing problem is to construct a gene combination (column) of size between $\ubar{k}$ and $\bar{k}$ with positive reduced cost that can be added to the restricted master problem or to prove that no such gene combination exists. Solving the pricing problem to optimality ensures that all gene combinations with positive reduced cost are identified.

Let $\pi_t\geq 0$, $\mu_n\geq 0$ and $\lambda\geq 0$ denote the dual variables corresponding to constraints (\ref{form2:cons1}), (\ref{form2:cons2}) and (\ref{form2:cons3}), respectively. The reduced cost of a gene combination $c$ is
\begin{align*}
   RC(c) = \sum_{t\in T_c} \pi_t - \sum_{n\in N_c} \mu_n - \lambda,
\end{align*}
where $T_c$ and $N_c$ are tumor and normal samples covered by gene combination $c$.

Let $u_g$ be a binary decision variable that takes the value 1 if gene $g$ is selected. Furthermore, we introduce binary decision variables $v_t$ and $v_n$ to indicate whether tumor sample $t$ or normal sample $n$ is covered by the selected genes. Recall that $\mathbf{M}_{sg}=1$ if sample $s$ has a mutated gene $g$. The pricing problem can be formulated as the following MIP
\begin{align}
    \max_{}\ & \sum_{t\in T} \pi_t v_t - \sum_{n\in N} \mu_n v_n - \lambda ,\label{formpp:cons1} \\
    \text{s.t.}\ &   \sum_{g\in G} u_g \geq \ubar{k}, \label{formpp:cons2} \\
    &  \sum_{g\in G} u_g \leq \bar{k} ,\label{formpp:cons3} \\
    & v_t \leq 1 - \frac{1}{\bar{k}} \sum_{g\in G: \mathbf{M}_{tg}=0} u_g, && \forall t\in T,\label{formpp:cons4} \\
    & v_n \geq 1 - \sum_{g\in G: \mathbf{M}_{ng}=0} u_g , && \forall n\in N,\label{formpp:cons6} \\
    & v_t \in \{0,1\} ,&&\forall t\in T,\label{formpp:cons7}  \\
    & v_n \in \{0,1\} ,&&\forall n\in N,\label{formpp:cons8}  \\
    & u_g \in \{0,1\} ,&&\forall g\in G.\label{formpp:cons9} 
\end{align}
The objective (\ref{formpp:cons1}) maximizes the reduced cost of the constructed gene combination. Constraints (\ref{formpp:cons2})-(\ref{formpp:cons3}) ensure that the constructed gene combination contains between $\ubar{k}$ and $\bar{k}$ genes. Constraints (\ref{formpp:cons4}) ensure that a tumor sample is covered only if it contains all genes in the constructed gene combination. If any gene is missing, the right-hand side becomes less than 1, forcing the binary variable to 0, while the constant $\frac{1}{\bar{k}}$ ensures the right-hand side remains non-negative. Constraints (\ref{formpp:cons6}) enforce that when a normal sample is covered by the constructed gene combination, the corresponding variable $v_n$ must take value 1, while it is lower bounded by 0 otherwise. A similar constraint is not required for tumor samples, due to the objective function. Finally, the domain of the variables is enforced by (\ref{formpp:cons7})-(\ref{formpp:cons9}).

\subsubsection{Solving the Pricing Problem}
We propose to solve the pricing problem using two methods. First, we solve it heuristically using a variable neighborhood search (VNS) (see \Cref{alg:VNS} in the Appendix). If the heuristic cannot identify any positive reduced cost combination, then we solve the pricing problem using an MIP solver. To accelerate the MIP solver, we propose the following changes. First, we return a solution as soon as a positive reduced cost combination is identified. Second, we observe that certain genes are frequently selected across successive pricing problem solutions. We fix $u_g=0$ for genes that were never selected in previous solutions. If this restricted pricing problem fails to produce a gene combination with positive reduced cost, we remove the restriction and solve the full pricing problem, ensuring that optimality is preserved. These steps are summarized in \Cref{alg:PP}.

Let $G_{r}\subseteq G$ be a restricted gene set, containing all genes that were in a generated column with positive reduced cost before. Initially, $G_r=\emptyset$. Furthermore, we restrict the search of the variable neighborhood search to a set $G_{\text{VNS}}\subseteq G$ containing the top $\gamma_3=1000$ genes. First, we apply a VNS heuristic on gene set $G_{\text{VNS}}\cup G_{r}$. Second, if that does not result in a positive reduced cost combination, we run the more expensive MIP solver on the restricted gene set $G_r$. Third, if no gene combination with positive reduced cost has been identified we apply the MIP solver on the full gene set $G$. Since the last step solves the pricing problem in an exact way, this procedure will either return a positive reduced cost combination or prove that no such combination exists. 

\begin{algorithm}
\caption{Exact pricing problem procedure}
\label{alg:PP}
\begin{algorithmic}[1]
\Require Gene sets $G_{r}$ and $G_{\text{VNS}}$
\Ensure A gene combination $c$ with positive reduced cost

\State $c\gets$ apply variable neighborhood search restricted to $G_{\text{VNS}}\cup G_{r}$ (see \Cref{alg:VNS})
\Statex \hspace{\algorithmicindent} \textbf{if} $c \neq \emptyset$ \textbf{goto} step 4.
\State $c\gets$ apply an MIP solver restricted to $G_{r}$
\Statex \hspace{\algorithmicindent} \textbf{if} $c \neq \emptyset$ \textbf{goto} step 4.
\State $c\gets$ apply an MIP solver using entire gene set $G$
\State \Return $c$
\end{algorithmic}
\end{algorithm}
\FloatBarrier

\subsection{Summary of the Column Generation-Based Framework}

\Cref{alg:colgen} summarizes the proposed column generation-based framework and shows how the greedy procedure, rounding heuristic, and pricing problem can be combined. The algorithm returns a set of selected gene combinations together with a valid lower bound ($LB^*$) and upper bound ($UB^*$) on the objective value. 

When only the initialization and integer version of the restricted master problem are performed (lines 1, 2 and 9), the algorithm reduces to the GR-CP and GR-MIP heuristics, depending on the solver employed. Conversely, executing all steps except the greedy procedure (line 2), yields our price-and-branch heuristic. This framework allows individual components to be enabled or added, supporting balancing computational cost and solution quality. 

\begin{algorithm}
\caption{Column generation-based framework}
\label{alg:colgen}
\begin{algorithmic}[1]
\Require Binary matrix $\mathbf{M}$
\Ensure Selected gene combinations, and bounds $LB^*$ and $UB^*$

\State initialize the restricted master problem (\ref{form2:cons1})-(\ref{form2:cons7}), with $\tilde{C}=\emptyset$

\State $\tilde{C}$ $\gets$ \texttt{GenerateGeneCombinations} (see \Cref{alg:greedyprocedure})

\Repeat
    \State $UB$ $\gets$ solve the restricted master problem
    \State $LB^*$ $\gets$ $\max\{LB^*, \texttt{RoundingHeuristic}\}$ (see \Cref{alg:roundingheuristic})
    \State $\tilde{C}$ $\gets$ \texttt{PricingProblemProcedure} (see \Cref{alg:PP})
\Until{no gene combination has positive reduced cost}

\State $UB^* \gets UB$

\State $LB_{} \gets$ solve the restricted master problem, with binary variables $z_c, \forall c\in \tilde{C}$

\State $LB^* \gets \max\{LB^*, LB_{}\}$

\end{algorithmic}
\end{algorithm}
\FloatBarrier

\begin{table}[h]
\caption{Summary of Datasets A and B. For each instance, we report the number of genes, samples, and mutations, as well as the required hit ranges.}
\label{tab:summarystatistics}
\centering
\begin{tabular}{lllrrr}
\toprule
dataset & cancer type & \#hits        & \#genes & \#samples & \#mutations \\ \midrule
\multirow{2}{*}{A}& BRCA        & 3, 2–3, 6–7 & 19,411 & 1242    & 837,207    \\
& PRAD        & 6–7, 7      & 17,274 & 752     & 476,115    \\ 
\midrule
\multirow{16}{*}{B}& BLCA-B        & 4, 6–7, 7   & 19,408 & 705     & 503,615    \\
& BRCA-B        & 3, 2–3, 6–7 & 20,728 & 1260    & 969,254    \\
& CESC-B        & 5, 6–7      & 20,073 & 613     & 430,358    \\
& COAD-B        & 3, 4, 5, 6  & 20,940 & 730     & 737,797    \\
& GBM-B        & 2, 3        & 19,828 & 693     & 931,859    \\
& HNSC-B        & 4, 5        & 19,573 & 806     & 656,209    \\
& KIRP-B        & 8           & 19,409 & 559     & 373,052    \\
& LGG-B         & 3, 6        & 19,167 & 821     & 468,223    \\
& LIHC-B        & 8           & 19,600 & 644     & 387,015    \\
& LUAD-B        & 3, 6–7      & 19,723 & 759     & 658,650    \\
& LUSC-B        & 5           & 19,672 & 647     & 467,554    \\
& OV-B          & 2, 6–7      & 19,984 & 663     & 682,506    \\
& PRAD-B        & 6–7, 7      & 18,978 & 770     & 446,045    \\
& SARC-B        & 5, 8        & 19,581 & 550     & 316,986    \\
& THCA-B        & 5, 6        & 19,305 & 756     & 505,205    \\
& UCEC-B        & 2, 6–7      & 21,200 & 835     & 1,021,354  \\ 
\bottomrule
\end{tabular}
\end{table}
\FloatBarrier

\section{Computational Results}\label{section4:results}
We evaluate our methods for the MHCDSCP on real-world cancer genomics data, demonstrating significantly improved computational efficiency, comparable accuracy to the state-of-the-art, and novel insights into solution quality.

\subsection{Experimental Setup}

\subsubsection{Characteristics of the Real-World Datasets}\label{section:exp_datasets}
Our experiments are conducted on two real-world somatic mutation datasets, Dataset A and~B. Originally extracted from The Cancer Genome Atlas (TCGA) \cite{weinstein2013cancer}, both datasets were provided by the authors of \cite{oles2025bigpicc} who processed the data according to the methodology in \cite{dash2019differentiating}. The characteristics of each dataset are summarized in \Cref{tab:summarystatistics} in terms of the numbers of genes, samples, and mutations. Both Dataset A and B contain 331 normal samples. The instances in Dataset B are distinguished by the ``-B'' suffix. Hit ranges are user-specified and are the same as those reported in \cite{oles2025bigpicc}. For example, a hit range of 2-3 indicates that a selected gene combination contains between 2 and 3 genes.

Note that Dataset A is publicly available via the Oak Ridge National Laboratory Gitlab\footnote{\url{https://code.ornl.gov/vo0/bigpicc}} and includes two cancer types that exactly match those used in~\cite{oles2025bigpicc}. Dataset B, including 16 cancer types, was provided by the authors of~\cite{oles2025bigpicc} in transposed matrix form and despite consultation could not be exactly matched. As a result, to ensure fairness, we use Dataset A for the direct comparison to the current state-of-the-art (\Cref{section:exp_a}), while Dataset B is used to demonstrate application to a wider range of cancer types (\Cref{section:exp_b}). For transparency, our processing steps to construct the input matrix for Dataset B are described in Appendix~\ref{appendix:A}.

\subsubsection{Technical Specifications \& Parameters}
All experiments are performed on a single quadcore Intel Core i7 CPU at 1.8 GHz and 16GB RAM. The CP formulation is solved using the open source solver OR-Tools 9.12, while all MIP formulations are solved using the commercial solver CPLEX 20.10. Both solvers use all available cores by default. We have also publicly released our code\footnote{\url{https://github.com/rickwillemsen/MHCDSCP}}.

Based on initial experiments, we fix the following parameters for all computational experiments. The number of selected gene combinations is limited to $\beta=10$ (preliminary results show that increasing $\beta$ does not have a large impact on the results). For both GR-CP and GR-MIP, we select the $\gamma_1=100$ genes that are most frequently mutated in tumor samples and generate at most $\gamma_2 = 100,000$ gene combinations. Price-and-branch is initialized with an empty set of columns, $\tilde{C}=\emptyset$.

The proposed methods may obtain an optimal integral solution, corresponding to an optimal binary classification with respect to the objective function. However, in practice, we are interested in how well a solution generalizes to unseen samples. To evaluate this, we partition the input data into a training (75\%) and a test (25\%) set, similar to \cite{oles2025bigpicc}. The mathematical models are trained on the training set and performance measures are computed on the test set. We evaluate classification performance using specificity, sensitivity, F1 score, and the Matthews Correlation Coefficient (MCC) (see Appendix~\ref{section2.3:performancemeasures} for formulas). The MCC is a symmetric measure and serves as the primary evaluation criterion in the current state-of-the-art method \cite{oles2025bigpicc}. Note that test set performance may differ even when an optimal solution is found for the training set. 

\subsubsection{Reporting of the Computational Results}
The performance measures in computational result tables are abbreviated as MCC, spec, sens, and F1. Additionally, we denote the number of selected gene combinations and corresponding objective value as \#comb and obj, respectively.

Solving the linear relaxation to optimality using price-and-branch provides an upper bound on the optimal objective value, which we use to compute an optimality gap for our methods (labeled bound and gap in tables, respectively) as follows:
\begin{align*}
    \text{gap} = \frac{\text{upper bound} - \text{objective value}}{|\text{objective value}|} \times 100.
\end{align*}

Computational time limits are established based on preliminary experiments to achieve comparable performance across models. We set a time limit of 20 minutes for GR-CP, 30~seconds for GR-MIP, and 5 minutes for price-and-branch. These limits exclude the time to initialize the models with greedy combinations for GR-CP and GR-MIP, which takes less than 30 seconds for all instances. In our tables, we report either the solver time (excluding initialization) or the total time (including initialization). Notably, no specific time limit is imposed on the MIP solver for the pricing problem.

\subsection{Case Study on Dataset A: Publicly Available Data}\label{section:exp_a}

First, we compare the performance of our three proposed heuristics: GR-CP, GR-MIP and price-and-branch. Second, we compare GR-MIP with the current state-of-the-art method. 

\subsubsection{Comparison Between GR-CP, GR-MIP, and Price-and-Branch}
\Cref{tab:CP} reports the results of GR-CP on Dataset A. While GR-CP reaches the time limit of 20 minutes for all five instances, it yields a feasible solution in every case. Note that the current GR-CP formulation outperformed several alternative CP formulations in initial experiments.

\begin{table}[h]
\caption{Performance of GR-CP on Dataset A.}
\label{tab:CP}
\centering
\resizebox{\textwidth}{!}{%
\begin{tabular}{llrrrrrrrrrrrrrr}
\toprule
\shortstack[l]{cancer\\ type} & \shortstack[l]{hit\\ range} &  & \multicolumn{3}{c}{solution} &  & \multicolumn{4}{c}{training performance} &  & \multicolumn{4}{c}{test performance} \\  \cline{4-6} \cline{8-11} \cline{13-16} 
 &          &  & \#comb     & obj     & \shortstack[l]{\tabtop solver\\ time (s)}   &  & MCC      & spec      & sens     & F1     &  & MCC     & spec    & sens    & F1     \\ \midrule
\multirow{3}{*}{BRCA}  & 2-3 &  & 10 & 677 & 1200.0 &  & 0.984 & 0.996 & 0.993 & 0.996 &  & 0.872 & 0.940 & 0.952 & 0.964 \\
     & 3   &  & 10 & 672 & 1200.0 &  & 0.973 & 0.984 & 0.991 & 0.993 &  & 0.887 & 0.940 & 0.961 & 0.969 \\
     & 6-7 &  & 10 & 505 & 1200.0 &  & 0.625 & 0.915 & 0.779 & 0.861 &  & 0.558 & 0.855 & 0.763 & 0.841 \\ \midrule
\multirow{2}{*}{PRAD}  & 6-7 &  & 10 & 239 & 1200.0 &  & 0.752 & 0.972 & 0.782 & 0.867 &  & 0.596 & 0.952 & 0.629 & 0.754 \\
     & 7   &  & 10 & 211 & 1200.0 &  & 0.667 & 0.952 & 0.712 & 0.814 &  & 0.525 & 0.952 & 0.543 & 0.687 \\ \midrule
avg  &     &  &  10.0  &     & 1200.0 &  & 0.800 & 0.964 & 0.851 & 0.906 &  & 0.688 & 0.928 & 0.769 & 0.843 \\ \bottomrule
\end{tabular}
}
\end{table}

The results of the same experiment using GR-MIP are presented in \Cref{tab:MIP}. GR-MIP achieves better objective values for all instances compared to GR-CP. Additionally, it successfully terminates within the time limit for two instances, leading to an average computation time of 26 seconds. This is faster than GR-CP, suggesting that an MIP model is preferred for this specific problem. However, when extending the problem with additional real-world (possibly nonlinear) constraints, a CP formulation may offer more flexibility and performance than an MIP formulation.

Some training performance measures reach their maximum value of 1 (e.g. spec). While this indicates that the obtained solution correctly classifies all training samples, the subsequent decline in performance on the test set highlights over-fitting, suggesting a training-test split may not be suitable for obtaining generalizable combinations. However, existing methods lack the efficiency to support more complex procedures. To ensure a fair comparison to prior work, we also adopt a training-test split. In \Cref{section5:conclusion} we identify potential directions enabled by our approach to mitigate overfitting.

\begin{table}[h]
\caption{Performance of GR-MIP on Dataset A.}
\label{tab:MIP}
\centering
\resizebox{\textwidth}{!}{%
\begin{tabular}{llrrrrrrrrrrrrrr}
\toprule
\shortstack[l]{cancer\\ type} & \shortstack[l]{hit\\ range} &  & \multicolumn{3}{c}{solution} &  & \multicolumn{4}{c}{training performance} &  & \multicolumn{4}{c}{test performance} \\\cline{4-6} \cline{8-11} \cline{13-16} 
&         &  & \#comb     & obj     & \shortstack[l]{\tabtop solver\\ time (s)}   &  & MCC      & spec      & sens     & F1     &  & MCC     & spec    & sens    & F1     \\ \midrule
\multirow{3}{*}{BRCA} & 2-3 &  & 10   & 678 & 30.0 &  & 0.986 & 1.000 & 0.993 & 0.996 &  & 0.878 & 0.928 & 0.961 & 0.967 \\
     & 3   &  & 10   & 677 & 30.0 &  & 0.984 & 1.000 & 0.991 & 0.996 &  & 0.864 & 0.928 & 0.952 & 0.962 \\
     & 6-7 &  & 10   & 517 & 30.0 &  & 0.663 & 0.968 & 0.773 & 0.866 &  & 0.551 & 0.880 & 0.737 & 0.828 \\ \midrule
\multirow{2}{*}{PRAD} & 6-7 &  & 10   & 244 & 17.3 &  & 0.775 & 0.988 & 0.788 & 0.877 &  & 0.642 & 0.964 & 0.667 & 0.787 \\
     & 7   &  & 10   & 216 & 22.5 &  & 0.692 & 0.988 & 0.693 & 0.814 &  & 0.564 & 0.952 & 0.590 & 0.725 \\ \midrule
avg  &     &  & 10.0 &     & 26.0 &  & 0.820 & 0.989 & 0.848 & 0.910 &  & 0.700 & 0.930 & 0.781 & 0.854 \\
    \bottomrule
\end{tabular}
}
\end{table}
\FloatBarrier

\Cref{tab:colgen_datasetA} reports the results when applying price-and-branch. This approach yields better objective values for four instances compared to GR-MIP. Furthermore, price-and-branch solves the root node to optimality within the time limit of 5 minutes for two instances, leading to upper bounds that can be used to calculate an optimality gap. 

The results in \Cref{tab:colgen_datasetA} incorporate the pricing problem speedups described in \Cref{alg:PP}, namely applying a variable neighborhood search and solving the pricing problem on a restricted set. The effect of the speedups is reported in \Cref{tab:pricingproblem_speedups} in the Appendix, showing that the number of completed column generation iterations increases from an average of 12 to 127 when both speedups are applied. Furthermore, across all configurations and instances, the total time spent solving the restricted master problem remained below 0.3 seconds, while the remaining time is spent on the pricing problem.

While price-and-branch obtains better solutions, the test performance of all three heuristics remains comparable. Since GR-MIP is much faster, namely terminating in under one minute when including initialization, it is the preferred method for further comparisons against the current state-of-the-art. 

\begin{table}[h]
\caption{Performance of Price-and-Branch on Dataset A.}
\label{tab:colgen_datasetA}
\resizebox{\textwidth}{!}{%
\begin{tabular}{llrrrrrrrrrrrrrrr}
\toprule
\shortstack[l]{cancer\\ type} & \shortstack[l]{hit\\ range} &  & \multicolumn{4}{c}{solution} &  & \multicolumn{4}{c}{training performance} &  & \multicolumn{4}{c}{test performance} \\ \cline{4-7} \cline{9-12} \cline{14-17} 
            &           &  & \#comb & obj & bound & \shortstack[l]{\tabtop solver\\ time (s)} &  & MCC   & spec  & sens  & F1    &  & MCC    & spec  & sens  & F1    \\ \midrule
\multirow{3}{*}{BRCA}        & 2-3       &  & 10   & 677 & 683   & 51.8        &  & 0.984 & 1.000 & 0.991 & 0.996 &  & 0.896  & 0.952 & 0.961 & 0.971 \\
            & 3         &  & 10   & 678 & 683   & 56.3        &  & 0.986 & 0.996 & 0.994 & 0.996 &  & 0.843  & 0.880 & 0.961 & 0.958 \\
            & 6-7       &  & 10   & 637 &   -    & 300.0       &  & 0.891 & 0.972 & 0.949 & 0.969 &  & 0.796  & 0.952 & 0.895 & 0.936 \\ \midrule
\multirow{2}{*}{PRAD}        & 6-7       &  & 10   & 280 &   -    & 300.0       &  & 0.876 & 0.984 & 0.899 & 0.940 &  & 0.639  & 0.916 & 0.724 & 0.809 \\
            & 7         &  & 10   & 254 &    -   & 300.0       &  & 0.783 & 0.931 & 0.858 & 0.897 &  & 0.488  & 0.855 & 0.629 & 0.721 \\ \midrule
avg         &           &  & 10.0 &     &       & 201.6       &  & 0.904 & 0.977 & 0.938 & 0.960 &  & 0.732  & 0.911 & 0.834 & 0.879 \\ \bottomrule
\end{tabular}
}
\end{table}
\FloatBarrier

\subsubsection{Comparison to the State-of-the-Art}
In this section, we compare GR-MIP to BiGPICC \cite{oles2025bigpicc}, the current state-of-the-art graph-based search heuristic that was executed on 280 nodes (estimated 11,760 CPU cores) on the Summit supercomputer. In comparison, GR-MIP is run on a single commodity CPU within one minute. Recall that we use a training set of 75\% of the data with the remaining 25\% for the test set, similar to the setup for BiGPICC~\cite{oles2025bigpicc}. 

In \Cref{tab:heuristic_B}, we report the performance values of BiGPICC as reported in the original paper~\cite{oles2025bigpicc}. When comparing the number of selected gene combinations with GR-MIP, we observe that we reduce the average number of gene combinations from 31 to 10. Despite the reduction in the number of gene combinations, performance measures remain similar. For instance, the average MCC is 0.694 and 0.700, respectively. Minor variations can be attributed to the differences in training-test split.

Although a direct comparison with BiGPICC is limited to Dataset A, it is clear that our GR-MIP heuristic is extremely fast, while maintaining a similar classification performance.

\begin{table}[h]
\caption{Comparison of test set performance between BiGPICC (results reported from \cite{oles2025bigpicc}) and GR-MIP on Dataset A.}
\label{tab:heuristic_B}
\centering
\resizebox{\textwidth}{!}{%
\begin{tabular}{lllrrrrrrrrrrr}
\toprule
\shortstack[l]{cancer\\ type} & \shortstack[l]{hit\\ range} &  & \multicolumn{5}{c}{\shortstack[l]{BiGPICC on Dataset A\\ (results reported from \cite{oles2025bigpicc})}}                                 &  & \multicolumn{5}{c}{GR-MIP on Dataset A}                   \\ \cline{4-8} \cline{10-14} 
            &           &  & \#comb & MCC   & spec & sens & F1    &  & \#comb & MCC   & spec & sens & F1    \\ \midrule
\multirow{3}{*}{BRCA}        & 2-3       &  & 5               & 0.900 & 0.892       & 0.987       & 0.974 &  & 10              & 0.878 & 0.928       & 0.961       & 0.967 \\
            & 3         &  & 9               & 0.900 & 0.892       & 0.987       & 0.974 &  & 10              & 0.864 & 0.928       & 0.952       & 0.962 \\
            & 6-7       &  & 34              & 0.749 & 0.916       & 0.882       & 0.922 &  & 10              & 0.551 & 0.880       & 0.737       & 0.828 \\ \midrule
\multirow{2}{*}{PRAD}        & 6-7       &  & 45              & 0.436 & 0.831       & 0.600       & 0.692 &  & 10              & 0.642 & 0.964       & 0.667       & 0.787 \\
            & 7         &  & 62              & 0.487 & 0.928       & 0.533       & 0.671 &  & 10              & 0.564 & 0.952       & 0.590       & 0.725 \\ \midrule
avg         &           &  & 31.0            & 0.694 & 0.892       & 0.798       & 0.847 &  & 10.0            & 0.700 & 0.930       & 0.781       & 0.854 \\ \bottomrule
\end{tabular}
}
\end{table}
\FloatBarrier

\begin{table}[H]
\caption{Test set performance of GR-MIP and Price-and-Branch on Dataset B.}
\label{tab:heuristic_A}
\centering
\resizebox{\textwidth}{!}{%
\begin{tabular}{llrrrrrrrrrrrrrr}
\toprule
\shortstack[l]{cancer\\ type} & \shortstack[l]{hit\\ range} &  & \multicolumn{8}{c}{GR-MIP on Dataset B}                                               &  & \multicolumn{4}{c}{Price-and-Branch on Dataset B}                \\ \cline{4-11} \cline{13-16} 
            &           &  & obj & gap   & \shortstack[l]{\tabtop total\\ time (s)} & \#comb & MCC   & spec  & sens  & F1    &  & obj & bound & gap  &  \shortstack[l]{total\\ time (s)}  \\ \midrule
\multirow{3}{*}{BLCA-B}        & 4         &  & 280       & 0.00  & 18.9       & 10   & 0.932 & 0.952 & 0.979 & 0.968 &  & 280       & 280   & 0.00 & 42.7       \\
        & 6-7       &  & 280       & 0.00  & 39.0       & 10   & 0.831 & 0.988 & 0.840 & 0.908 &  & 279       & 280   & 0.36 & 100.9      \\
        & 7         &  & 270       &   -    & 39.4       & 10   & 0.863 & 1.000 & 0.862 & 0.926 &  & 276       &    -   &   -   & 311.7      \\ \midrule
\multirow{3}{*}{BRCA-B}        & 2-3       &  & 697       & 0.00  & 52.4       & 9    & 0.984 & 0.988 & 0.996 & 0.996 &  & 697       & 697   & 0.00 & 46.0       \\
        & 3         &  & 697       & 0.00  & 47.2       & 7    & 0.967 & 0.964 & 0.996 & 0.991 &  & 697       & 697   & 0.00 & 49.8       \\
        & 6-7       &  & 645       &    -   & 45.9       & 10   & 0.728 & 0.988 & 0.819 & 0.898 &  & 688       &   -    &   -   & 300.7      \\ \midrule
\multirow{2}{*}{CESC-B}        & 5         &  & 211       & 0.00  & 14.9       & 10   & 0.974 & 1.000 & 0.972 & 0.986 &  & 211       & 211   & 0.00 & 15.3       \\
        & 6-7       &  & 211       & 0.00  & 33.9       & 10   & 0.840 & 1.000 & 0.817 & 0.899 &  & 211       & 211   & 0.00 & 14.4       \\ \midrule
\multirow{4}{*}{COAD-B}        & 3         &  & 296       & 1.01  & 50.8       & 10   & 0.934 & 0.964 & 0.970 & 0.970 &  & 299       & 299   & 0.00 & 69.1       \\
        & 4         &  & 296       & 1.01  & 59.2       & 10   & 0.935 & 0.988 & 0.950 & 0.969 &  & 299       & 299   & 0.00 & 64.3       \\
        & 5         &  & 286       & 4.55  & 29.9       & 10   & 0.848 & 1.000 & 0.850 & 0.919 &  & 298       & 299   & 0.34 & 191.9      \\
        & 6         &  & 284       & 5.28  & 28.3       & 10   & 0.886 & 1.000 & 0.890 & 0.942 &  & 298       & 299   & 0.34 & 213.1      \\ \midrule
\multirow{2}{*}{GBM-B}         & 2         &  & 271       & 0.00  & 1.5        & 5    & 0.989 & 1.000 & 0.989 & 0.994 &  & 271       & 271   & 0.00 & 11.6       \\
         & 3         &  & 271       & 0.00  & 29.9       & 9    & 0.955 & 0.952 & 1.000 & 0.978 &  & 271       & 271   & 0.00 & 17.5       \\ \midrule
\multirow{2}{*}{HNSC-B}        & 4         &  & 356       & 0.00  & 20.7       & 10   & 0.928 & 0.940 & 0.983 & 0.971 &  & 356       & 356   & 0.00 & 47.7       \\
       & 5         &  & 356       & 0.00  & 24.0       & 10   & 0.920 & 0.976 & 0.950 & 0.966 &  & 356       & 356   & 0.00 & 44.5       \\ \midrule
KIRP-B        & 8         &  & 159       & 7.55  & 36.1       & 10   & 0.763 & 1.000 & 0.702 & 0.825 &  & 170       & 171   & 0.59 & 120.4      \\ \midrule
\multirow{2}{*}{LGG-B}         & 3         &  & 367       & 0.00  & 22.9       & 10   & 0.929 & 0.940 & 0.984 & 0.972 &  & 367       & 367   & 0.00 & 14.7       \\
         & 6         &  & 357       & 2.80  & 41.2       & 10   & 0.893 & 0.976 & 0.927 & 0.954 &  & 366       & 367   & 0.27 & 63.3       \\ \midrule
LIHC-B        & 8         &  & 224       & 4.91  & 38.6       & 10   & 0.774 & 1.000 & 0.744 & 0.853 &  & 234       & 235   & 0.43 & 225.2      \\ \midrule
\multirow{2}{*}{LUAD-B}        & 3         &  & 321       & 0.00  & 23.6       & 10   & 0.936 & 0.940 & 0.991 & 0.972 &  & 321       & 321   & 0.00 & 25.9       \\
        & 6-7       &  & 321       & 0.00  & 18.3       & 10   & 0.907 & 0.988 & 0.925 & 0.957 &  & 321       & 321   & 0.00 & 58.0       \\ \midrule
LUSC-B        & 5         &  & 237       & 0.00  & 12.4       & 10   & 0.988 & 0.988 & 1.000 & 0.994 &  & 237       & 237   & 0.00 & 15.1       \\ \midrule
\multirow{2}{*}{OV-B}          & 2         &  & 248       & 0.40  & 1.5        & 9    & 0.904 & 0.952 & 0.952 & 0.952 &  & 249       & 249   & 0.00 & 11.4       \\
          & 6-7       &  & 244       & 2.05  & 43.2       & 10   & 0.964 & 1.000 & 0.964 & 0.982 &  & 248       & 249   & 0.40 & 186.4      \\ \midrule
\multirow{2}{*}{PRAD-B}        & 6-7       &  & 325       & 1.23  & 37.9       & 10   & 0.822 & 0.976 & 0.855 & 0.913 &  & 329       & 329   & 0.00 & 71.1       \\
        & 7         &  & 317       & 3.79  & 40.0       & 10   & 0.804 & 1.000 & 0.809 & 0.894 &  & 329       & 329   & 0.00 & 108.6      \\ \midrule
\multirow{2}{*}{SARC-B}        & 5         &  & 164       & 0.00  & 5.1        & 10   & 0.926 & 1.000 & 0.909 & 0.952 &  & 164       & 164   & 0.00 & 31.8       \\
        & 8         &  & 146       & 12.33 & 21.0       & 10   & 0.785 & 1.000 & 0.727 & 0.842 &  & 163       & 164   & 0.61 & 113.6      \\ \midrule
\multirow{2}{*}{THCA-B}        & 5         &  & 319       & 0.00  & 19.0       & 10   & 0.909 & 1.000 & 0.915 & 0.956 &  & 319       & 319   & 0.00 & 31.7       \\
        & 6         &  & 319       & 0.00  & 36.6       & 10   & 0.897 & 0.988 & 0.915 & 0.951 &  & 319       & 319   & 0.00 & 86.7       \\ \midrule
\multirow{2}{*}{UCEC-B}        & 2         &  & 378       & 0.00  & 1.2        & 10   & 0.970 & 0.988 & 0.984 & 0.988 &  & 378       & 378   & 0.00 & 29.3       \\
        & 6-7       &  & 362       & 4.42  & 41.6       & 10   & 0.887 & 0.988 & 0.913 & 0.950 &  & 378       & 378   & 0.00 & 87.9       \\ \midrule
avg         &           &  &           & 1.66  & 29.6       & 9.7  & 0.896 & 0.982 & 0.911 & 0.945 &  &           &       & 0.11 & 85.5   \\ \bottomrule
\end{tabular}
}
\end{table}
\FloatBarrier

\subsection{Case Study on Dataset B: Solution Quality and Scalability}\label{section:exp_b}
We investigate the performance of GR-MIP on additional cancer types. We focus on Dataset~B, containing 16 cancer types from TCGA. The results are presented in \Cref{tab:heuristic_A}. GR-MIP achieves high values of MCC for all instances within a minute. On average, we achieve an MCC of 0.896. This performance is better compared to the MCC values reported for BiGPICC (see Table 2 in \cite{oles2025bigpicc}), which achieved an average MCC value of 0.687 when run on supercomputing infrastructure. However, a direct comparison is limited since our Dataset B could not be exactly matched to datasets used in \cite{oles2025bigpicc} (as detailed in~\Cref{section:exp_datasets}). 

Furthermore, our price-and-branch solves the root node to optimality within 5 minutes for all but two instances. The resulting valid upper bounds allow for the calculation of optimality gaps for both GR-MIP and price-and-branch. Across the 33 benchmark instances, GR-MIP solves 18 instances to optimality, whereas price-and-branch solves 23 instances to provable optimality. 

These results indicate that it is possible to solve instances of the MHCDSCP to provable optimality. Additionally, our GR-MIP heuristic is able to find near-optimal solutions to real-world instances within one minute.

\section{Conclusion}\label{section5:conclusion}
Existing approaches to identify carcinogenic multi-hit gene combinations rely on exhaustive search-based (but still heuristic) methods and require supercomputing infrastructure. In contrast, we introduce fast and practical methods that obtain comparable solutions on a single commodity CPU in minutes. Specifically, we introduce a CP and MIP formulation for the problem. Additionally, we propose a price-and-branch heuristic, which allows us to solve instances to provable optimality and demonstrate that our MIP heuristic provides near-optimal solutions within a minute. These results suggest that identifying these gene combinations is far less computationally intensive than previously believed.

Given that the problem is now computationally tractable to solve, this opens up the possibility for other research directions. It becomes possible to test hypotheses in multi-hit theory that were previously intractable. Since our proposed methods are flexible they can easily accommodate a range of practical assumptions, such as changing the weight of tumor and normal samples, requiring all tumor samples to be covered, varying the number of selected gene combinations, or fixing the size of each gene combination (hit range). These optimization models are easily accessible to researchers in cancer genomics.

For future work, a potential direction is designing an exact method guaranteed to solve the problem to optimality. The key bottleneck is the reliance on a MIP solver to prove optimality of the pricing problem; a dedicated algorithm could accelerate this bottleneck. Additionally, investigating complex cross-validation procedures to mitigate any occurrences of overfitting may be warranted, especially given that this is now possible due to the significantly improved computational performance of our approach. While this work prioritizes computational efficiency and optimality insights for the existing modeling of the problem (motivated by the significant drawbacks of prior parallelization-reliant approaches), the CP formulation provides an extensible foundation for future research to explore complex biological constraints that may be challenging to express in a linearized formulation.

%%
%% Bibliography
%%

%% Please use bibtex, 

\bibliography{lipics-v2021-sample-article}

\appendix

\section{Performance Measures}\label{section2.3:performancemeasures}
To evaluate classification performance, we use several measures, namely specificity, sensitivity, F1 score, and the Matthews Correlation Coefficient (MCC)~\cite{chicco2020advantages}. These measures are defined as
\begin{align*}
    \text{sensitivity} &= \frac{\text{TP}}{\text{TP} + \text{FN}},\\
    \text{specificity} &= \frac{\text{TN}}{\text{TN} + \text{FP}},\\
    \text{precision} &= \frac{\text{TP}}{\text{TP} + \text{FP}},\\
    \text{F1 score} &= 2\cdot \frac{\text{precision} \cdot \text{sensitivity}}{\text{precision} + \text{sensitivity}},\\
    \text{MCC} &= \frac{\text{TP} \cdot \text{TN} - \text{FP} \cdot \text{FN}}{\sqrt{ (\text{TP} + \text{FP})(\text{TP} + \text{FN})(\text{TN} + \text{FP})(\text{TN} + \text{FN})}}.
\end{align*}
Note that precision is not used as a performance measure, but is used to define the F1 score. 

\section{Algorithmic Details for the Pricing Problem}\label{appendix:PP}

Variable neighborhood search (VNS) can be used to quickly identify solutions to the pricing problem and is summarized in \Cref{alg:VNS}. A classic VNS consists of three steps:
\begin{enumerate}
    \item \textbf{Shake:} Move to a random solution in the $\kappa$-th neighborhood by randomly replacing $\kappa$~genes in the current combination with genes from $G$ to escape local optima.
    \item \textbf{Local search:} A series of systematic swaps of each gene in the current combination $c$ with a gene in a set $g\in G\setminus c$. The swap is accepted, only if it improves the reduced cost. This step is used to find a local optimum.
    \item \textbf{Neighborhood change:} If no improvement is found, the search space is expanded to a larger neighborhood ($\kappa \gets \kappa + 1$). When an improvement is found, the first neighborhood is used again ($\kappa=1$).
\end{enumerate}
To increase computational efficiency, we restrict the search to a set $G_{\text{VNS}}\subseteq G$ containing the top $\gamma_3$ genes, where we set $\gamma_3=1000$. We terminate the classic VNS after $\iota=5$ iterations. Since in our problem, the combination should satisfy a hit size between $\ubar{k}$ and $\bar{k}$, we run the classic VNS for each possible hit size $k$.

\begin{algorithm}[H]
\caption{Variable neighborhood search}
\label{alg:VNS}
\begin{algorithmic}[1]
\Require A gene set $G'$
\Ensure A set of gene combinations $c^*$ with positive reduced cost

\State $C^*\gets \emptyset$
\For{$k\in [\ubar{k}, \bar{k} ] \cap \mathbb{Z}$}
\State $c\gets$ generate a random gene combination of size $k$
\State $i\gets 1$
\While{$i \leq \iota$}
\State $\kappa\gets 1$
\While{$\kappa \leq k$}
\State $c'\gets \text{Shake}(c, \kappa, G')$: Randomly replace $\kappa$ genes by genes in $G'$.
\State $c''\gets \text{LocalSearch}(c', G')$: Replace each gene by genes in $G'$.
\If{$RC(c'')> RC(c)$}
\State $c\gets c''$
\State $\kappa \gets 1$
\Else
\State $\kappa \gets \kappa + 1$
\EndIf
\EndWhile
\State $i\gets i+1$
\EndWhile
\If{$RC(c)>0$}
\State $C^*\gets C^* \cup c$
\EndIf
\EndFor
\end{algorithmic}
\end{algorithm}
\FloatBarrier

\section{Construction of Dataset B}\label{appendix:A}
As discussed in~\Cref{section:exp_datasets}, Dataset B was provided by the authors of \cite{oles2025bigpicc} in a partially processed transposed form. According to \cite{oles2025bigpicc}, the provided dataset was originally obtained from TCGA and processed according to the methodology described in~\cite{dash2019differentiating} to obtain only the somatic mutations in tumor and normal samples. To transform the provided datasets to binary input matrix form for our Dataset B, we perform the following steps:

\begin{itemize}

\item The somatic mutations in normal samples are provided as a list of ($g$, $s$)-tuples, indicating that gene $g$ is mutated in sample $s$. This can be converted into a binary matrix and requires no further processing.

\item The partially processed tumor data is also provided as a list of ($g$, $s$, $i$)-tuples, where $i$ denotes the number of experiments in which gene $g$ was observed to be mutated in sample~$s$. We consider a gene to be mutated if at least one experiment reports a mutation, i.e. $i\geq 1$ and set the corresponding entry in the binary matrix to $\mathbf{M}_{sg}=1$. Any gene~$g$ that is never mutated in the tumor data and does not appear in the normal data is removed, as it does not contribute to the problem.

\end{itemize}

\section{Accelerating the Pricing Problem}
\Cref{tab:pricingproblem_speedups} reports the performance of price-and-branch with and without speedups mentioned in \Cref{alg:PP}. Without any speedups, price-and-branch manages to solve on average 12 column generation iterations. When first solving a reduced MIP formulation, the number of iterations increases on average to 104, indicating that the pricing problem is solved faster. When adding the final speedup, which is a variable neighborhood search (VNS), the number of iterations increases further to 127 on average. Furthermore, this is the only configuration that solves the root node to optimality for some instances. In most cases, activating a speedup also results in an improved objective value.
\begin{table}[h]
\caption{Effect of the pricing problem speedups (corresponding to the three steps in \Cref{alg:PP}) on Dataset A. We report the number of column generation iterations (\#iter). The solver time (in seconds) is split into its two components: solving the restricted master problem (MP) and solving the pricing problem (PP).}\label{tab:pricingproblem_speedups}
\resizebox{\textwidth}{!}{%
\begin{tabular}{llrrrrrrrrrrrrrrr}
\toprule
\shortstack[l]{cancer\\ type} & \shortstack[l]{hit\\ range} &  & \multicolumn{4}{c}{full MIP}            &  & \multicolumn{4}{c}{reduced MIP + full MIP}             &  & \multicolumn{4}{c}{VNS + reduced MIP + full MIP}         \\ \cline{4-7} \cline{9-12} \cline{14-17} 
            &           &  & obj & \#iter & \shortstack[l]{\tabtop MP\\ time (s)}  &\shortstack[l]{PP\\ time (s)} &  & obj & \#iter & \shortstack[l]{MP\\ time (s)} & \shortstack[l]{PP\\ time (s)} &  & obj & \#iter &\shortstack[l]{MP\\ time (s)}& \shortstack[l]{PP\\ time (s)} \\ \midrule
\multirow{3}{*}{BRCA}        & 2-3       &  & 665 & 8          & 0.0    & 300.0  &  & 669       & 15         & 0.0    & 300.0  &  & 677       & 62         & 0.1    & 51.5 \\
            & 3         &  & 651 & 8          & 0.0    & 300.0  &  & 661       & 24         & 0.0    & 300.0  &  & 678       & 113        & 0.1    & 56.1 \\
            & 6-7       &  & 590 & 10         & 0.0    & 300.0  &  & 615       & 87         & 0.1    & 299.9  &  & 637       & 115        & 0.2    & 299.8  \\ \midrule
\multirow{2}{*}{PRAD}        & 6-7       &  & 264 & 17         & 0.0    & 300.0  &  & 281       & 177        & 0.2    & 299.8  &  & 280       & 212        & 0.3    & 299.7  \\
            & 7         &  & 235 & 16         & 0.0    & 300.0  &  & 257       & 217        & 0.2    & 299.8  &  & 254       & 134        & 0.1    & 299.9  \\ \midrule
avg         &           &  &     & 11.8       & 0.0    & 300.0  &  &           & 104.0      & 0.1    & 299.9  &  &           & 127.2      & 0.2    & 201.4  \\
\bottomrule
\end{tabular}
}
\end{table}
\FloatBarrier

\end{document}